\documentclass{article}
 \usepackage{amsmath}

\newtheorem{thm}{Theorem}[section]
\newtheorem{cor}[thm]{Corollary}

\newtheorem{lem}[thm]{Lemma}

\newtheorem{rem}[thm]{Remark}

\newtheorem{ex}[thm]{Example}

\newcommand{\be}{\begin{equation}}
\newcommand{\ee}{\end{equation}}
\newcommand{\ben}{\begin{enumerate}}
\newcommand{\een}{\end{enumerate}}
\newcommand{\beq}{\begin{eqnarray}}
\newcommand{\eeq}{\end{eqnarray}}
\newcommand{\beqn}{\begin{eqnarray*}}
\newcommand{\eeqn}{\end{eqnarray*}}

\newcommand{\pa}{\partial}

\newcommand{\pxi}{ {\pa \over \pa x^i}}

\newcommand{\pyj}{ {\pa \over \pa y^j}}

\newcommand{\qed}{\hspace*{\fill}Q.E.D.}  

\begin{document}
\title{Conformal Vector Fields of a Class of  Finsler Spaces II}
\author{Guojun Yang\footnote{Supported by the
National Natural Science Foundation of China (11471226) }}
\date{}
\maketitle
\begin{abstract}

In this paper, we first give two fundamental principles under a
technique to characterize conformal vector fields of
$(\alpha,\beta)$ spaces to be homothetic and determine the local
structure of those homothetic fields. Then we use the principles
to study conformal vector fields of some classes of
$(\alpha,\beta)$ spaces under certain curvature conditions.
Besides, we construct a family of non-homothetic conformal vector
fields on a family of locally projectively Randers spaces.


 {\bf MR(2000) subject classification: }
53B40
\end{abstract}

\section{Introduction}

Let $F$ be a Finsler metric on a manifold $M$ with the fundamental
metric matrix $(g_{ij})$, and $V$ be a vector field on $M$. $V$ is
called a conformal vector field of the Finsler manifold $(M,F)$ if
$V=V^i\pa/\pa x^i$ satisfies
 \be\label{y01}
 V_{0|0}=-2cF^2,
 \ee
where the symbol $_|$ denotes the horizontal covariant derivative
with respect to Cartan, or Berwald or Chern connection, and
$V_{0|0}=V_{i|j}y^iy^j$, $V_i=g_{im}V^m$ and $c$, called a
conformal factor, is a scalar function on $M$. If $c=constant$,
then $V$ is called homothetic. If $c=0$, $V$ is called a Killing
vector field. An equivalent description of (\ref{y01}) is
(\ref{j8}) below (\cite{HM1}).

Conformal vector fields (esp. Killing vector fields) play an
important role in Finsler geometry. When $F$ is a Riemann metric,
it is shown that the local solutions of a conformal vector field
can be determined on $(M,F)$ if $(M,F)$ is of constant sectional
curvature in dimension $n\ge 3$ (\cite{YShen} \cite{ShX}), or more
generally locally conformally flat in dimension $n\ge 2$
(\cite{Y1}), or under other curvature conditions (\cite{Y2}). Some
important problems in Finsler geometry can be solved by
constructing a conformal vector field of a Riemann metric with
certain curvature features (especially of constant sectional
curvature) (cf. \cite{BRS} \cite{Y3}--\cite{Y9} \cite{Yu1}).

A vector field on a manifold $M$ induces a flow $\varphi_t$ acting
on $M$, and $\varphi_t$ is naturally lifted to a flow
$\widetilde{\varphi}_t$ on the tangent bundle $TM$, where
$\widetilde{\varphi}_t: TM \mapsto TM$ is defined  by
  $\widetilde{\varphi}_t(x,y)=(\varphi_t(x),\varphi_{t*}(y))$, where $x\in M, y\in T_xM$. In
\cite{MH}, Huang-Mo define a homothetic vector field on a Finsler
space by
 \be\label{y002}
 \widetilde{\varphi}_t^*F=e^{-2ct}F,
 \ee
where $c$ is a constant. Then in \cite{MH}, it obtains the
relation between the flag curvatures of two Finsler metric $F$ and
$\tilde{F}$, where $\tilde{F}$ is defined by $(F,V)$ under
navigation technique for a homothetic vector field $V$ of $F$.
Note that for a scalar function $c$, (\ref{y01}) does not imply
(\ref{y002}). We show in \cite{Y8} that for a vector field $V$ and
a scalar function $c$ on a Finsler manifold $(M,F)$ , if
(\ref{y002}) holds, then $c(\varphi_t(x))=c(x)$; and (\ref{y01})
implies (\ref{y002}) iff. $c(\varphi_t(x))=c(x)$.

 An {\it $(\alpha,\beta)$-metric} is defined by a
Riemannian metric
 $\alpha=\sqrt{a_{ij}(x)y^iy^j}$ and a $1$-form $\beta=b_i(x)y^i$ on a manifold
 $M$, which can be expressed in the following form:
 $$F=\alpha \phi(s),\ \ s=\beta/\alpha,$$
where $\phi(s)$ is a function satisfying certain conditions such
that $F$ is regular. In \cite{SX}, Shen-Xia study conformal vector
fields of Randers spaces under certain curvature conditions. In
\cite{K}, Kang characterizes conformal vector fields of
$(\alpha,\beta)$-spaces by some PDEs in a special case
$\phi'(0)\ne 0$. Later on, we prove the same result for all
non-Riemannian $(\alpha,\beta)$-spaces (without the condition
$\phi'(0)\ne 0$) (\cite{Y8}). In this paper, we mainly show some
curvature conditions of some $(\alpha,\beta)$-spaces on which
every conformal vector field must be homothetic, and also consider
the possible local solutions for those homothetc vector fields.

First, we show two fundamental results on conformal vector fields
to be homothetic and the local solutions for homothetc vector
fields respectively.

\begin{thm}\label{th1}
 Let $F=\alpha\phi(\beta/\alpha)$ be an
  $(\alpha,\beta)$-metric on an  $n$-dimensional
 manifold $M$, where $\phi(s)$ satisfies $\phi(s)\ne \sqrt{1+ks^2}$ for any constant $k$ and
 $\phi(0)=1$. Define a
 Riemann metric $h=\sqrt{h_{ij}(x)y^iy^j}$ and a 1-form $\rho=p_i(x)y^i$ by
  \be\label{y1}
h=\sqrt{u(b^2)\alpha^2+v(b^2)\beta^2},\ \ \rho=w(b^2)\beta,\ \
 (b^2:=||\beta||_{\alpha}^2),
  \ee
 where $u=u(t)\ne 0,v=v(t),w=w(t)\ne 0$ are some smooth
 functions. Suppose $\rho$ is a conformal 1-form of $h$. Then any
 conformal vector field of $(M,F)$ must be homothetic.
\end{thm}

\begin{thm}\label{th2}
 In Theorem \ref{th1}, further  suppose the dimension $n\ge 3$ and $h$ is of constant sectional curvature
 $\mu$. Locally we can express
  \be\label{y2}
 h=\frac{2}{1+\mu|x|^2}|y|, \ \ \ \ \
 \rho=\frac{4}{(1+\mu|x|^2)^2}\Big\{-2\big(\lambda+\langle
 d,x\rangle\big)\langle x,y\rangle+|x|^2\langle d,y\rangle+p_r^ix^ry^i+\langle e,y\rangle\Big\},
 \ee
 where $\lambda$ is a constant number, $d,e$ are constant vectors and
  $P=(p_i^j)$ is a skew-symmetric matrix. Let $V=V^i(x)\pa/\pa x^i$ be a conformal vector field of
  $(M,F)$ with the conformal factor $c$. Then we have one of the following cases:
   \ben
  \item[{\rm (i)}] ($\mu=0,\lambda=0,d=0$)  $V$ is given by
   \beq
 &&V^i=-2\tau x^i+q^i_mx^m+\gamma^i,\label{y3}\\
 &&\hspace{-0.4cm}\big(Qe=P\gamma,\ PQ-QP=2\tau P\big).\label{y03}
   \eeq
   In this case,   we have  $c=\tau$.

\item[{\rm (ii)}] ($\mu=0,\lambda\ne0$) $V$ is given by
  \beq
&& V^i=q^i_mx^m+\gamma^i,\label{y4} \\
  &&\hspace{-2.8cm}\big(\langle d,\gamma\rangle=0,\ Qd=0,\ Qe=-2\lambda\gamma+P\gamma,\ PQ-QP=2R\big),\label{y04}
   \eeq
   where $R=(r^i_j)$ is defined by $r^i_j=\gamma^id^j-\gamma^jd^i$. In this case,   we have
  $c=0$, and so $V$ is a
 Killing vector field.

\item[{\rm (iii)}] ($\mu\ne 0$) $V$ is given by
 \beq
&& V^i=2\mu\langle\gamma,x\rangle
 x^i+(1-\mu|x|^2)\gamma^i+q^i_mx^m,\label{y5}\\
&&\hspace{-2.8cm}\big(\langle d+\mu e,\gamma\rangle=0,\
Qd=-\mu(2\lambda\gamma+P\gamma),\ Qe=-2\lambda\gamma+P\gamma,\
PQ-QP=2R\big),\label{y05}
 \eeq
 where $R=(r^i_j)$ is defined by
 $r^i_j=\mu(e^i\gamma^j-e^j\gamma^i)+\gamma^id^j-\gamma^jd^i$. In this case,   we have  $c=0$, and so $V$ is a
 Killing vector field.
 \een
 In (\ref{y3})--(\ref{y05}), $\tau$ is a constant number, $\gamma=(\gamma^i)$ is a constant
 vector
 and  $Q=(q^i_k)$ is a constant skew-symmetric matrix.
\end{thm}

 Note that for $\lambda\ne
0$ in (\ref{y05}), we have
$$Qd=-\mu(2\lambda\gamma+P\gamma),\ Qe=-2\lambda\gamma+P\gamma \ \ \Longrightarrow \ \   \langle d+\mu e,\gamma\rangle=0.$$

In Theorem \ref{th2}, if additionally $\rho$ is closed, then $V$
is given by (\ref{y3})--(\ref{y05}) with $P=0$ and $d=\mu e$. This
case has actually been obtained in \cite{Y8}. See Corollary
\ref{cor41} below.

In Theorem \ref{th2}, if  $\rho$ is a homothetic 1-form of $h$,
then $V$ is given by (\ref{y3})--(\ref{y05}) with $d=0$ and
$\mu=0$, or $d=-\mu e$ and $\lambda=0$. See Corollary \ref{cor42}
below.

In Theorem \ref{th2}, if the condition that $\rho$ is a conformal
1-form of $h$ is replaced by the condition that $\rho$ is closed,
then a conformal vector field of $(M,F)$ is not necessarily
homothetic (see Remark \ref{rem73}). This implies that for a
locally projectively flat Randers space $(M, F)$ with
$F=\alpha+\beta$, a conformal vector field of $(M,F)$ is not
necessarily homothetic (cf. \cite{Y8}). We will show a detailed
construction of such a family of examples in Section \ref{sec7}.

\

Under some certain curvature conditions of an
$(\alpha,\beta)$-space, usually we can define $h$ and $\rho$ by
(\ref{y1}) by choosing suitable functions $u,v$ and $w$, such that
$\rho$ is a conformal 1-form of $h$, or $\rho$ is a conformal
1-form of $h$ and $h$ is of constant sectional curvature.
Following this idea, in Sections \ref{sec5}--\ref{sec7}, as an
application of Theorem \ref{th1} and Theorem \ref{th2}, we will
study some properties of conformal vector fields of some special
$(\alpha,\beta)$ spaces under certain curvature conditions. In
Section 5, we study conformal vector fields on
$(\alpha,\beta)$-spaces of isotropic S-curvature, and our main
result is Theorem \ref{th51}. In Section 6, we study conformal
vector fields on $(\alpha,\beta)$-spaces of Douglas type, and our
main result is Theorem \ref{th61}. In Section 7, we study
conformal vector fields on locally projectively flat Randers
spaces, and our main result is Example \ref{ex72}.

\section{Preliminaries}

Let $F$ be a Finsler metric on  a manifold $M$, and $V$ be a
vector field on $M$. Let $\varphi_t$ be the flow generated by $V$.
Define $\widetilde{\varphi}_t: TM \mapsto TM$ by
  $\widetilde{\varphi}_t(x,y)=(\varphi_t(x),\varphi_{t*}(y))$.
For a conformal vector field $V$ defined by (\ref{y01}), Huang-Mo
show in  \cite{HM1} an equivalent definition in the way
 \be\label{j8}
 \widetilde{\varphi}_t^*F=e^{-2\sigma_t}F,
 \ee
where $\sigma_t$ is a function on $M$ for every $t$, and in this
case, $c$ in (\ref{y01})and $\sigma_t$ in (\ref{j8}) are related
by
 $$
  \sigma_t=\int_0^t c(\varphi_s)ds, \ \ \ c=\frac{d}{dt}_{|t=0}\sigma_t.
 $$
 Now in a special case of (\ref{j8}), suppose
$\widetilde{\varphi}_t^*F=e^{-2ct}F$, namely,
 \be\label{y5}
 F(\varphi_t(x),\varphi_{t*}(y))=e^{-2c(x)t}F(x,y),
 \ee
 where $c$ is a scalar function on $M$.
 Differentiating (\ref{y5}) by $t$ at $t=0$,  we obtain
  \be\label{y6}
 X_V(F)=-2cF,
  \ee
  where
   \be\label{XV}
   X_V:=V^i\pxi+y^i\frac{\pa
   V^j}{\pa x^i}\pyj
   \ee
   is a vector field on $TM$. However, (\ref{y6}) generally does
   not imply (\ref{y5}) for a scalar function (see \cite{Y8})). We have a different description for conformal
   vector fields.

\begin{lem} (\cite{Y8})\label{lem20}
A vector
   field $V$ on a Finsler manifold $(M,F)$ is conformal iff.
    \be\label{g14}
  X_V(F^2)=-4cF^2  \ (\Longleftrightarrow
 X_V(F)=-2cF),
   \ee
   where $c$ is a scalar function on $M$.
\end{lem}

For conformal vector fields on a Riemann manifold, we have the
following lemma.

\begin{lem} (\cite{Y1} \cite{Y8})\label{lem21}
Let $\alpha$ be a Riemann metric of constant sectional curvature
$\mu$ on an $n$-dimensional manifold $M$. Locally express $\alpha$
by
 $$
\alpha=\frac{2}{1+\mu|x|^2}|y|.
 $$
 \ben
  \item[{\rm (i)}] ($n\ge 3$) Let $V$ be a
conformal vector field of $(M,\alpha)$ with the conformal factor
$c=c(x)$. Then locally we have
 \beqn
 V^i=-2\big(\lambda+\langle
 d,x\rangle\big)x^i+|x|^2d^i+q_r^ix^r+\eta^i,\ \
c=\frac{\lambda(1-\mu |x|^2)+ \langle \mu \eta+d,x\rangle}{1+\mu
 |x|^2},
 \eeqn
where $\lambda$ is a constant number, $d,\eta$ are constant
vectors and
  $(q_i^j)$ is skew-symmetric.

\item[{\rm (ii)}] ($n\ge 2$)  In (i), if additionally the 1-form
$V_iy^i$ is closed ($V_i:=a_{ik}V^k$), then locally we have
 $$
 V^i=-2(\lambda+\mu \langle
 e,x\rangle) x^i+(1+\mu|x|^2)
 e^i,\ \ \ c=\frac{\lambda(1-\mu |x|^2)+ 2\mu\langle e,x\rangle}{1+\mu
 |x|^2}.
 $$
 \een
\end{lem}

For an $(\alpha,\beta)$ space $(M,F)$ with
$F=\alpha\phi(\beta/\alpha)$, a conformal vector field is
characterized by the following lemma, which is also proved in
\cite{K} by assuming $\phi'(0)\ne 0$.

\begin{lem}(\cite{Y8})\label{lem22}
  Let $F=\alpha\phi(\beta/\alpha)$ be an
  $(\alpha,\beta)$-metric on an  $n$-dimensional
 manifold $M$, where $\phi(s)$ satisfies $\phi(s)\ne \sqrt{1+ks^2}$ for any constant $k$ and
 $\phi(0)=1$.
 Then $V=V^i(x)\pa/\pa x^i$ is a conformal vector field of $(M,F)$ with the
conformal factor $c=c(x)$ if and only if
  \be\label{y13}
  V_{i;j}+V_{j:i}=-4ca_{ij},\ \ \ \ V^jb_{i;j}+b^jV_{j;i}=-2cb_i,
  \ee
where $V_i$ and $b^i$ are defined by $V_i:=a_{ij}V^j$ and
$b^i:=a^{ij}b_j$, and the covariant derivatives are taken with
respect to the Levi-Civita connection of $\alpha$.

\end{lem}

 In this paper, for a Riemannian metric $\alpha
=\sqrt{a_{ij}y^iy^j}$ and a $1$-form $\beta = b_i y^i $,  let
 $$r_{ij}:=\frac{1}{2}(b_{i|j}+b_{j|i}),\ \ s_{ij}:=\frac{1}{2}(b_{i|j}-b_{j|i}), \ \ s_j:=b^is_{ij}, $$
 where we define $b^i:=a^{ij}b_j$, $(a^{ij})$ is the inverse of
 $(a_{ij})$, and $\nabla \beta = b_{i|j} y^i dx^j$  denotes the covariant
derivatives of $\beta$ with respect to $\alpha$.

\section{Proof of Theorem \ref{th1}}

To prove Theorem \ref{th1}, we first show two lemmas in Riemann
geometry.

\begin{lem}\label{lem31}
The Lie bracket of two conformal vector fields of a Riemann
manifold is a conformal vector.
\end{lem}

\begin{lem}\label{lem32}
Let $V$ and $cV$ be two conformal vector fields of a Riemann
manifold, where $c=c(x)$ is a scalar function on the manifold.
Then $c$ must be a constant.
\end{lem}

In \cite{HM1}, Huang-Mo have  generalized Lemmas \ref{lem31} and
\ref{lem32}  to general Finsler manifolds.

For a pair $(\alpha,\beta)$, define a new pair $(h,\rho)$ by
(\ref{y1}), we obtain in \cite{Y8}, an equivalent characterization
of (\ref{y13}) in terms of $(h,\rho$) as follows.

\begin{lem}\cite{Y8}\label{lem33}
Let $\alpha=\sqrt{a_{ij}y^iy^j}$ be a Riemann metric and
 $\beta=b_iy^i$ be a 1-form and $V=V^i\pa/\pa x^i$ be a vector field
  on an $n$-dimensional manifold $M$.  Define a
 Riemann metric $h=\sqrt{h_{ij}(x)y^iy^j}$ and a 1-form $\rho=p_i(x)y^i$
 by (\ref{y1}), namely,
  $$
h=\sqrt{u(b^2)\alpha^2+v(b^2)\beta^2},\ \ \rho=w(b^2)\beta,\ \
 (b^2:=||\beta||_{\alpha}^2),
  $$
 where $u=u(t)\ne 0,v=v(t),w=w(t)\ne 0$ are some smooth functions.
  Then $\alpha$, $\beta$  and $V$ satisfy (\ref{y13}) if and
 only if
  \be\label{y14}
 V_{0|0}=-2ch^2,\ \ \ \ V^jp_{i|j}+p^jV_{j|i}=-2cp_i,
  \ee
  where $p^j:=h^{ij}p_i,\ V_j:=h_{ij}V^i$ and the covariant derivative is taken
  with respect to the Levvi-Civita connection of $h$.
\end{lem}

Now we show the proof of Theorem \ref{th1}. Since $V$ is a
conformal vector field of the $(\alpha,\beta)$ space $(M,F)$ with
$F=\alpha\phi(\beta/\alpha)$, by Lemma \ref{lem22}, $\alpha$,
$\beta$  and $V$ satisfy (\ref{y13}) with the conformal factor
$c=c(x)$. Then by Lemma \ref{lem33}, $h$, $\rho$  and $V$ satisfy
(\ref{y14}). The first formula in (\ref{y14}) is equivalent to
$V_{i|j}+V_{j|i}=-4ch_{ij}$, and  then using this, a simple
observation shows that the second formula in (\ref{y14}) is
equivalent to
 $$[V,W]=2cW,$$
where $W=W^i\pa/\pa x^i$ is a vector field defined by
$W^i:=h^{ij}p_j$. The first  formula in (\ref{y14}) shows that $V$
is a conformal vector field of $h$, and by assumption, $W$ is also
a conformal vector field of $h$. Now by Lemma \ref{lem31}, $[V,W]$
is a conformal vector field of $h$, and so is $cW$. Thus by Lemma
\ref{lem32}, we have $c=constant$, which means that $V$ is a
homothetc vector field of $(M,F)$.         \qed

\section{Proof of Theorem \ref{th2}}

In Theorem \ref{th2}, since $\rho$ is a conformal 1-form of $h$
and  $V$ is a conformal vector field of
  $(M,F)$ with the conformal factor $c$, we have $c=constant$ by
  Theorem \ref{th1}. Actually, we can also directly verify
  $c=constant$ in giving the local solution of $V$ under the assumption of Theorem \ref{th2}.

Since $\rho$ is a conformal 1-form of $h$ and  $V$ is a conformal
vector field of  $(M,F)$, we have (\ref{y14}) by Lemma \ref{lem22}
and Lemma \ref{lem33}. By  the assumption that the dimension $n\ge
3$ and $h$ is of  constant sectional curvature $\mu$, it follows
from Lemma \ref{lem21} (i) that locally we can write
 \be\label{y15}
 p_i=\frac{4}{(1+\mu|x|^2)^2}\Big\{-2\big(\lambda+\langle
 d,x\rangle\big)x^i+|x|^2d^i+p_r^ix^r+e^i\Big\},
 \ee
and correspondingly  $V$ and $c$ are  expressed as
 \beq
 V^i\hspace{-0.6cm}&&=-2\big(\tau+\langle
 \eta,x\rangle\big)x^i+|x|^2\eta^i+q_r^ix^r+\gamma^i,\label{y16}\\
c\hspace{-0.6cm}&&=\frac{\tau(1-\mu |x|^2)+ \langle \mu
\gamma+\eta,x\rangle}{1+\mu |x|^2},\label{y17}
 \eeq
where $\lambda,\tau$ are a constant numbers, $d,e,\eta,\gamma$ are
constant vectors and $P=(p^i_j)$,
  $Q=(q_i^j)$ are  skew-symmetric constant matrices.
The second equation in (\ref{y14}) is equivalent to
 \be\label{y18}
 V^j\frac{\pa p_i}{\pa
 x^j}+p_j\frac{\pa V^j}{\pa x^i}=-  2cp_i.
 \ee
Now plugging (\ref{y15}), (\ref{y16}) and (\ref{y17}) into
(\ref{y18}) yields an equivalent equation of a polynomial in
$(x^i)$ of degree four, in which every degree must be zero. Then
we respectively have (from degree zero to degree four)
 \beq
Qe\hspace{-0.5cm}&&=-2\lambda \gamma+P\gamma,\label{y19}\\
 0\hspace{-0.5cm}&&=2\big(\langle
 \eta,e\rangle-\langle d,\gamma\rangle+2\lambda\tau\big)x^i-2(\langle
 e,x\rangle\eta^i+\langle d,x\rangle\gamma^i+\mu\langle \gamma,x\rangle e^i-\langle \gamma,x\rangle d^i)\nonumber\\
 &&+(p^i_jq^j_k-q^i_jp^j_k-2\tau p^i_k)x^k,\label{y20}\\
 0\hspace{-0.5cm}&&=\big[2\lambda(\eta^i-\mu\gamma^i)-\mu
 q_k^ie^k+4\tau(\mu e^i-d^i)-q^i_kd^k+p^i_k\eta^k+\mu p^i_k\gamma^k\big]|x|^2\nonumber\\
 &&+2\big[2\lambda\mu\langle\gamma,x\rangle+4\tau\langle
 d,x\rangle-\langle d,Qx\rangle+\langle \eta,Px\rangle\big]x^i-2\langle\mu\gamma+\eta,x\rangle p^i_kx^k,\label{y21}\\
 0\hspace{-0.5cm}&&=\Big\{2\mu(\langle\eta, e\rangle-\langle d,\gamma\rangle
 -2\lambda\tau)x^i-2\mu\langle
 e,x\rangle\eta^i-2\mu\langle
 d,x\rangle\gamma^i\nonumber\\
 &&+2\langle\eta, x\rangle(\mu e^i-d^i)+\mu(p^i_jq^j_k-q^i_jp^j_k+2\tau p^i_k)x^k\Big\}|x|^2
 +4\langle
 d,x\rangle\langle\eta+\mu\gamma,x\rangle x^i,\label{y22}\\
 0\hspace{-0.5cm}&&=\mu\big[(2\lambda
\eta^i-q^i_kd^k+p^i_k\eta^k)|x|^2-2\langle
 2\lambda\eta- Qd+P\eta,x\rangle x^i\big].\label{y23}
 \eeq
We will determine $V$ from (\ref{y19})--(\ref{y23}). First we will
show that $\mu=0$ and $\eta=0$, or $\tau=0$ and $\eta=-\mu
\gamma$. Then by (\ref{y17}) we have $c=constant$ given by
$c=\tau$ or $c=0$.

\

\noindent{\bf Case I :} Assume $\mu=0$. We will prove that
$\eta=0$ from (\ref{y19})--(\ref{y23}). Let $\eta\ne0$ in the
following discussion for this case. We show a contradiction.

Plugging $\mu=0$ into (\ref{y22}) yields
 $$
 -2\langle \eta,x\rangle(d^i|x|^2-2\langle d,x\rangle x^i)=0,
 $$
which implies $d=0$. Plugging $\mu=0$ and $d=0$ into (\ref{y21})
yields
 \be\label{y24}
 (p^i_k\eta^k+2\lambda\eta^i)|x|^2-2(\langle P\eta,x\rangle
 x^i+\langle\eta,x\rangle p^i_kx^k)=0.
 \ee
Contracting (\ref{y24}) by $x^i$ yields
  $$|x|^2\big(\langle2\lambda \eta-P\eta,x\rangle\big)=0,$$
which implies $P\eta=2\lambda \eta$. Since $P$ is real and
skew-symmetric, its real characteristic roots must be zeros. Thus
it follows from $P\eta=2\lambda \eta$ that $\lambda=0$ and
$P\eta=0$. Now plugging $\lambda=0$ and $P\eta=0$ back into
(\ref{y24}) gives $\langle\eta,x\rangle p^i_kx^k=0$. So we get
$P=0$. Further, plugging $\mu=0$, $d=0$, $\lambda=0$ and $P=0$
into (\ref{y20}), and then contracting it by $x^i$ we get
 $$
 \langle\eta,e\rangle|x|^2-\langle\eta,x\rangle \langle
 e,x\rangle=0,
 $$
which shows that $e=0$. Now since  $d=0$, $\lambda=0$, $e=0$ and
$P=0$, we obtain $\beta=\rho=0$ from (\ref{y15}).  Thus
$F=\alpha\phi(\beta/\alpha)=\alpha$ is Riemannian. It is a
contradiction.

\

\noindent{\bf Case II :} Assume $\mu\ne 0$. We will prove that
$\eta=-\mu \gamma$ and $\tau=0$ from (\ref{y19})--(\ref{y23}).

{\bf Case IIa :} Suppose  $\eta\ne -\mu \gamma$. We will show a
contradiction.

By (\ref{y23}) we get
 \be\label{y25}
 Qd=2\lambda\eta +P\eta.
 \ee
Contracting (\ref{y20}) by $x^i$ we get
 $$
 \big(2\lambda\tau-\langle d,\gamma\rangle+\langle
 \eta,e\rangle\big)|x|^2-\langle e,x\rangle\langle \eta+\mu\gamma,x\rangle=0,
 $$
which implies $e=0$ (since $\eta\ne -\mu \gamma$). By (\ref{y22}),
it directly shows that $d=0$ (since $\eta\ne -\mu \gamma$). Now by
(\ref{y25}), we have $P\eta=-2\lambda\eta$ which implies
 \be\label{y26}
 \eta=0, \ {\rm or} \ \lambda=0 \ {\rm and} \ P\eta=0,
 \ee
  since $P$ is skew-symmetric.
Similarly, by (\ref{y19}), we have $P\gamma=2\lambda \gamma$ which
implies
 \be\label{y27}
 \gamma=0, \ {\rm or} \ \lambda=0 \ {\rm and} \ P\gamma=0.
 \ee
Since $\eta\ne -\mu \gamma$, by (\ref{y26}) and (\ref{y27}) we get
$\lambda=0,P\eta=0$ and $P\gamma=0$. By this fact and $d=0,e=0$,
we obtain $-2\langle\eta+\mu \gamma,x\rangle Px=0$ from
(\ref{y21}), which shows $P=0$. Then again we get a contradiction
since $F$ is Riemannian in this case.

{\bf Case IIb :} By the discussion in Case IIa, we have obtained
 $\eta= -\mu \gamma$. Suppose $\tau\ne 0$. We will get a
 contradiction.

Plug $\eta= -\mu \gamma$ into (\ref{y21}), we get
 $$
 \big(4\mu\tau e^i-\mu q^i_ke^k-4\mu\lambda\gamma^i-4\tau
 d^i-q^i_kd^k\big)|x|^2
 +2\langle 4\tau
 d+2\mu\lambda\gamma+
 Qd+\mu P\gamma,x\rangle x^i=0,
 $$
from which we obtain
 \be\label{y28}
 Qd=4\mu\tau e-\mu Qe-4\mu\lambda\gamma-4\tau
 d,
 \ee
 \be\label{y29}
 Qd=-4\tau d-2\mu\lambda\gamma-\mu P\gamma.
 \ee
Now plugging (\ref{y19}) into $(\ref{y28})-(\ref{y29})$ gives
$4\mu\tau e=0$, which shows $e=0$. Then by (\ref{y19}) we have
$P\gamma=2\lambda \gamma$, which implies $\gamma=0$, or
$\lambda=0$ and $P\gamma=0$. In either case, we have $Qd=-4\tau d$
from (\ref{y29}). Thus we obtain $d=0$ since $\tau\ne0$. Plugging
$d=e=0$ into (\ref{y20}) and contracting it by $x^i$ we get
$4\lambda\tau|x|^2=0$, which shows $\lambda=0$. Now we have proved
 \be\label{y30}
 \lambda=0,\ \ d=0,\ \ e=0.
 \ee
Plug (\ref{y30}) into (\ref{y20}) and (\ref{y22}), we respectively
obtain
 \be\label{y31}
 PQx-QPx-2\tau Px=0,
 \ee
\be\label{y32}
 \mu(PQx-QPx+2\tau Px)|x|^2=0.
 \ee
Then $(\ref{y32})/(\mu|x|^2)-(\ref{y31})$ yields $4\tau Px=0$,
from which we get $P=0$. So $F$ is Riemannian again by (\ref{y30})
and $P=0$. This is a contradiction.

\

Therefore, we have $\mu=0$ and $\eta=0$, or $\eta=-\mu \gamma$ and
$\tau=0$ by the discussion in Case I and Case II above. Then we
will solve $V$ from (\ref{y19})--(\ref{y23}) in the following.

\

\noindent{\bf Case A :} Assume $\mu=0$ and $\eta=0$. Plugging
$\mu=0$ and $\eta=0$ into (\ref{y20}) gives
 \be\label{y33}
 2(2\lambda\tau-\langle d,\gamma\rangle)x^i-2\tau
 p^i_kx^k+p^i_jq^j_kx^k-q^i_jp^j_kx^k-2\langle
 d,x\rangle\gamma^i+2\langle
 \gamma,x\rangle d^i=0.
 \ee
Contracting (\ref{y33}) by $x^i$ we get
 \be\label{y34}
 \langle d,\gamma\rangle=2\lambda\tau.
\ee Then by (\ref{y34}), (\ref{y33}) is equivalent to
 \be\label{y35}
 2\tau P=PQ-QP-2R, \ \ \big(R=(r^i_j) \ {\rm with} \
 r^i_j:=\gamma^id^j-\gamma^jd^i\big)
 \ee
 Plugging
$\mu=0$ and $\eta=0$ into (\ref{y21}) gives
 $$
 -(q^i_kd^k+4\tau d^i)|x|^2+2(4\tau\langle d,x\rangle-\langle
 d,Qx\rangle)x^i=0,
 $$
which is equivalent to
 \be\label{y36}
 Qd=-4\tau d.
 \ee
 Since $Q$ is skew-symmetric, (\ref{y36}) is equivalent to
  \be\label{y37}
 d=0,\ \ {\rm or} \ \tau=0 \ {\rm and}  \ Qd=0.
  \ee
For (\ref{y22}) and (\ref{y23}), they automatically hold since
$\mu=0$ and $\eta=0$. Then (\ref{y19})--(\ref{y23}) are equivalent
to (\ref{y19}), (\ref{y34}), (\ref{y35}) and (\ref{y37}), which
are broken into one of the following three cases:
 \be\label{y38}
 (\mu=0,\eta=0),\ \ d=0,\ \ \lambda=0,\ \ Qe=P\gamma,\ \ 2\tau P=PQ-QP,
 \ee
 \be\label{y39}
  (\mu=0,\eta=0),\ \ d=0,\ \ \tau=0,\ \ Qe=-2\lambda\gamma+P\gamma,\ \ PQ=QP,
 \ee
\be\label{y40}
 (\mu=0,\eta=0),\ \  \tau=0,\ \ Qd=0,\ \  \langle d,\gamma\rangle=0,\ \ Qe=-2\lambda\gamma+P\gamma,\ \
 PQ-QP=2R.
 \ee
We see (\ref{y39}) $\Rightarrow$ (\ref{y40}). So we have only two
cases (\ref{y38}) and (\ref{y40}), which are just Theorem
\ref{th2} (i) and (ii) respectively.

\

\noindent{\bf Case B :} Assume $\eta=-\mu\gamma$ and $\tau=0$.
First (\ref{y23}) $\Leftrightarrow$ (\ref{y25}), which, by
$\eta=-\mu\gamma$,  is written as
 \be\label{y41}
 Qd=-\mu(2\lambda\gamma+P\gamma).
 \ee
Plugging $\eta=-\mu\gamma$ and $\tau=0$ into (\ref{y20}) yields
 \beq
 0\hspace{-0.5cm}&&=-2(\langle
 d,\gamma\rangle+\mu\langle\gamma,e\rangle)x^i+(p^i_jq^j_k-q^i_jp^j_k)x^k\nonumber\\
 &&+2\langle\gamma,x\rangle
 d^i-2\langle d,x\rangle \gamma^i
 +2\mu(\langle e,x\rangle\gamma^i-\langle \gamma,x\rangle e^i).\label{y42}
 \eeq
Then contracting (\ref{y42}) by $x^i$ gives
 \be\label{y43}
 \langle d,\gamma\rangle=-\mu\langle\gamma,e\rangle.
 \ee
By (\ref{y43}), (\ref{y42}) is equivalent to
 \be\label{y44}
 PQ-QP=2R,\ \ \big(R=(r^i_j) \ {\rm with}\
 r^i_j:=\mu(e^i\gamma^j-e^j\gamma^i)+\gamma^id^j-\gamma^jd^i\big).
 \ee
Now it follows that (\ref{y21}) automatically holds from
$\eta=-\mu\gamma$, $\tau=0$, (\ref{y19}) and (\ref{y41});
(\ref{y22}) automatically holds from $\eta=-\mu\gamma$, $\tau=0$,
(\ref{y43}) and (\ref{y44}). So (\ref{y19})--(\ref{y23}) are
equivalent to (\ref{y19}), (\ref{y41}), (\ref{y43}) and
(\ref{y44}). This gives Theorem \ref{th2} (iii). \qed

\

 We consider two special cases of Theorem \ref{th2}: (ia) $\rho$
 is additionally closed ($\Leftrightarrow$ $P=0$ and $d=\mu e$ in (\ref{y2}));
 (ib) $\rho$ is homothetic ($\Leftrightarrow$ $d=0$ and
$\mu=0$, or $d=-\mu e$ and $\lambda=0$ in (\ref{y2})).
  Then we obtain the following two
corollaries respectively.

\begin{cor} \label{cor41} In Theorem \ref{th2}, additionally assume $\rho$ is a
closed 1-form. Then $V$ is determined by Theorem \ref{th2} (i)
with (\ref{y03}) being replaced by $Qe=0$, or Theorem \ref{th2}
(ii) with (\ref{y04}) being replaced by $Qe=-2\lambda \gamma$, or
Theorem \ref{th2} (iii) with (\ref{y05}) being replaced by
$\langle\gamma,e\rangle=0$ and $Qe=-2\lambda \gamma$.
\end{cor}

\begin{cor}  \label{cor42}  In Theorem \ref{th2}, additionally assume $\rho$ is homothetic.
 Then $V$ is determined by Theorem \ref{th2} (i), or Theorem \ref{th2}
(ii) with (\ref{y04}) being replaced by
 $$
 Qe=-2\lambda\gamma+P\gamma,\ PQ-QP=0,
 $$
 or Theorem \ref{th2} (iii) with (\ref{y05})
being replaced by
 $$
Qe=P\gamma,\ PQ-QP=2R,
 $$
 where $R=(r^i_j)$ is defined by $r^i_j:=2\mu(e^i\gamma^j-e^j\gamma^i)$.
\end{cor}

\section{Isotropic S-curvature}\label{sec5}

 In this section, we use Theorems \ref{th1} and \ref{th2} to study conformal
 vector fields on all $(\alpha,\beta)$-spaces of
 isotropic S-curvature.

The S-curvature is one of the most important non-Riemannian
quantities in Finsler geometry which was originally introduced for
the volume comparison theorem (\cite{shen2}).  The S-curvature is
said to be isotropic if there is a scalar function
$\theta=\theta(x)$ on $M$ such that
$${\bf S}=(n+1)\theta F.$$
If $\theta$ is a constant, then we call $F$ is of constant
S-curvature.

In this section, we mainly prove the following theorem.

 \begin{thm}\label{th51}
 Let $F=\alpha \phi(s)$, $s=\beta/\alpha$, be an $n(\ge 2)$-dimensional
 non-Riemannian $(\alpha,\beta)$-metric on the manifold $M$, where $\phi(0)=1$.
   Suppose  $F$ is  of isotropic
   S-curvature. Then any conformal vector field of $(M,F)$ is
 homothetic.
\end{thm}

In \cite{HM1}, Huang-Mo prove Theorem \ref{th51} when
$F=\alpha+\beta$ is a Randers metric.  But our proof here for a
Randers metric is quite different from that in \cite{HM1}. To
prove Theorem \ref{th51}, we first show the following lemma.

\begin{lem}\label{lem51}
  Let $F=\alpha \phi(s)$, $s=\beta/\alpha$, be an $n(\ge 2)$-dimensional
non-Riemannian  $(\alpha,\beta)$-metric on the manifold $M$, where
$\phi(0)=1$.
   Suppose  $F$ is  of isotropic
   S-curvature. Then we have only three classes:
 \ben
   \item[{\rm (i)}] (\cite{CS}) ($n\ge 2$)  $F$ is of Randers type in the equivalent form
   $F=\alpha+\beta$ satisfying
 \be\label{y45}
 r_{ij}=2\theta(a_{ij}-b_ib_j)-b_is_j-b_js_i,
 \ee
 where $\theta$ is a scalar function. In this case, the
 S-curvature is given by ${\bf S}=(n+1)\theta F$.

 \item[{\rm (ii)}] (\cite{CS1}) ($n\ge 2$) $\phi(s)$ is arbitrary and $r_{ij}=0,s_i=0$. In this
 case, ${\bf S}=0$.

\item[{\rm (iii)}] (\cite{CSY}) ($n= 2$) $\phi(s)$ is given by
 \be\label{y46}
      \phi(s)=\big\{(1+k_1s^2)(1+k_2s^2)\big\}^{\frac{1}{4}}e^{\int^s_0\theta(s)ds},
     \ee
 and $\beta$ satisfies
  \be\label{y47}
  r_{ij}=\frac{3k_1+k_2+4k_1k_2b^2}{4+(k_1+3k_2)b^2}(b_is_j+b_js_i),
  \ee
   where $\theta(s)$ is defined by
      \be\label{y48}
       \theta(s):=\frac{\pm\sqrt{k_2-k_1}}
      {2(1+k_1s^2)\sqrt{1+k_2s^2}},
      \ee
      and $k_1$ and $k_2$ are constants with $k_2>k_1$. In this
 case, ${\bf S}=0$.
 \een
\end{lem}

We should  note that  the condition $r_{ij}=0,s_i=0$ is a special
case of (\ref{y47}). In Lemma \ref{lem51}, the second class is
almost trivial, and if $F$ is not of Randers type, the third class
is essential and in this case, the norm $||\beta||_{\alpha}$ may
not be a constant. In \cite{Y9}, we further determine the local
structure of the third class in Lemma \ref{lem51}, and it is an
Einstein metric but generally not Ricci-flat.

\

{\it Proof of Theorem \ref{th51} :} Define a Riemann metric $h$
and and a 1-form $\rho$ by (\ref{y1}).  We will show by
(\ref{y45}) or (\ref{y47}), there are  functions $u,v,w$ such that
$\rho$ is a conformal 1-form of $h$.

\

 \noindent {\bf Case I:}
 Suppose that (\ref{y45}) holds. In this case, define in (\ref{y1})
 \be\label{y49}
 u(b^2)=k_2w(b^2),\ \
 v(b^2)=\Big(k_1-\frac{k_1+k_2}{b^2}\Big)w(b^2),
 \ee
where $k_1,k_2$ are constant and $w(b^2)$ is some function such
that $h$ is a Riemann metric. For $u,v,w$ defined by (\ref{y49}),
the equation (\ref{y45}) is equivalently transformed into
 \be\label{y50}
 \widetilde{r}_{ij}=-\frac{2\theta}{k_1}\Big(\frac{1}{1-b^2}+\frac{b^2w'(b^2)}{w(b^2)}\Big)h_{ij},
 \ee
where we define $\widetilde{r}_{ij}:=(p_{i|j}+p_{j|i})/2$ for the
covariant derivative of $\rho$ with respect to $h$. So by
(\ref{y50}), $\rho$ is a conformal 1-form with respect to $h$.
Therefore, by Theorem \ref{th1}, we get the proof of Theorem
\ref{th51} in this case. In particular,  if $k_1,k_2$ and $w(b^2)$
are taken as $k_1=1,k_2=-1$ and $w(b^2)=b^2-1$  in(\ref{y49}),
then we get the navigation data for a Randers metric.

\

\noindent {\bf Case II:} Suppose that (\ref{y47}) holds. We have
shown in \cite{Y9} that by choosing
 $$u(b^2)=1,\ \ v(b^2)=0,\ \
 w(b^2)=(1+k_1b^2)^{-\frac{3}{4}}(1+k_2b^2)^{-\frac{1}{4}},$$
the equation (\ref{y47}) is equivalent to
 \be\label{y51}
 \widetilde{r}_{ij}=0,
 \ee
where we define $\widetilde{r}_{ij}:=(p_{i|j}+p_{j|i})/2$ for the
covariant derivative of $\rho$ with respect to $h$. The equation
(\ref{y51}) shows that $\rho$ is a Killing form of $h$ (a special
case of conformality). Thus by Theorem \ref{th1}, we get the proof
of Theorem \ref{th51} in this case.  \qed

\

A Randers metric will be of isotropic S-curvature under some
curvature conditions, for example, (ia) $F$ is a weak Einstein
metric (\cite{SYi}); (ib) $F$ is Ricci-reversible (\cite{SY1});
(ic) $F$ is of Einstein-reversibility (\cite{Y4}). Then we have
the following corollary.

\begin{cor}
 On a manifold $M$, if a Randers metric $F=\alpha+\beta$ is weakly Einsteinian, or Ricci-reversible, or
 of Einstein-reversibility, then any conformal vector field of $(M,F)$ is
 homothetic.
\end{cor}

In \cite{SYi}, Shen-Yidirim show that if a Randers metric is of
weakly isotropic flag curvature in dimension $n\ge 3$, then
 $h$ is of constant sectional curvature and $\rho$ is a
conformal 1-form  by choosing the navigation data $u(b^2)=1-b^2,\
v(b^2)=b^2-1,\
 w(b^2)=b^2-1$ in (\ref{y1}). Thus we have the following
corollary.

\begin{cor}
 On a manifold $M$ of dimension $n\ge 3$, if a Randers metric $F=\alpha+\beta$ is of weakly isotropic flag curvature,
 then any conformal vector field $V$ of $(M,F)$ is homothetic, and
 locally $V$ can be determined by Theorem \ref{th2}.
\end{cor}

In\cite{BRS}, Bao-Robles-Shen prove that if a Randers metric is of
constant flag curvature, then $h$ is of constant sectional
curvature and $\rho$ is a homothetic 1-form by choosing the
navigation data $u(b^2)=1-b^2,\ v(b^2)=b^2-1,\
 w(b^2)=b^2-1$ in (\ref{y1}).   Thus
we have  the following corollary.

\begin{cor}
 On a manifold $M$ of dimension $n\ge 3$, if a Randers metric $F=\alpha+\beta$ is of constant flag curvature,
 then any conformal vector field $V$ of $(M,F)$ is homothetic, and
 locally $V$ can be determined by Corollary \ref{cor42}.
\end{cor}

\section{Douglas metrics}\label{sec6}

In this section, we study conformal vector fields of
$(\alpha,\beta)$-spaces of Douglas type. We obtain the following
theorem.

\begin{thm}\label{th61}
Let $F=\alpha \phi(s)$, $s=\beta/\alpha$, be an $n(\ge
2)$-dimensional
  $(\alpha,\beta)$-metric of non-Randers type on the manifold $M$, where $\phi(0)=1$.
   Suppose  $F$ is  Douglas metric. Then any conformal vector field of $(M,F)$ is
 homothetic.
\end{thm}

A Randers metric $F=\alpha+\beta$ is a Douglas metric iff. $\beta$
is closed (\cite{BaMa}). A conformal vector field on a Randers
space of Douglas type is not necessarily homothetic (see Example
\ref{ex72} below). To prove Theorem \ref{th61}, we first show the
following lemma.

\begin{lem}
\label{lem61}
 Let $F=\alpha \phi(s)$, $s=\beta/\alpha$, be an $n$-dimensional
  $(\alpha,\beta)$-metric, where $\phi(0)=1$. Suppose that
    $\beta$ is not parallel with respect to $\alpha$ and $F$ is not
of Randers type. Then we have
   \ben
  \item[{\rm (i)}] (\cite{LSS})  ($n\ge 3$) $F$ is  a Douglas metric if and only if
    \beq
    &&\big\{1+(k_1+k_3)s^2+k_2s^4\big\}\phi''(s)=(k_1+k_2s^2)\big\{\phi(s)-s\phi'(s)\big\},\nonumber\\
    && b_{i|j}=\theta
     \big\{(1+k_1b^2)a_{ij}+(k_2b^2+k_3)b_ib_j\big\},\label{y52}
   \eeq
 where  $k_1,k_2,k_3$ are constant with $k_2\ne k_1k_3$, and $\theta$ is a scalar function.

 \item[{\rm (i)}]  (\cite{Y3}) ($n=2$) $F$ is  a Douglas metric if and only if
  \beq
 && F \ {\rm is \ of \ the \ metric \ type} \
F_0^{\pm}=\alpha\pm \beta^2/\alpha,\nonumber\\
 && r_{ij}=\theta\big\{(1\pm 2b^2)a_{ij}\mp
  3b_ib_j\big\}
  +\frac{3}{\pm1-b^2}(b_is_j+b_js_i),\label{y53}
  \eeq
 where $\theta$ is a scalar function.
 \een
\end{lem}

{\it Proof of Theorem \ref{th61} :} Define a Riemann metric $h$
and and a 1-form $\rho$ by (\ref{y1}).  We will show by
(\ref{y52}) or (\ref{y53}), there are  functions $u,v,w$ such that
$\rho$ is a conformal 1-form of $h$.

\

 \noindent {\bf Case I:}
 Suppose that (\ref{y52}) holds. In this case, let $u,v,w$  in (\ref{y1})
be defined by
 \beqn
&&u(b^2)=1,\ \ v(b^2)=0, \\
 &&w(b^2)=e^{-\int^{b^2}_0\frac{1}{2}\frac{k_3+k_2t}{1+(k_1+k_3)t+k_2t^2}dt}
 \eeqn
Then it is shown in \cite{Y3}  that (\ref{y52}) is equivalent to
 \be\label{y54}
 p_{i|j}=\theta w(b^2)(1+k_1b^2)h_{ij}.
 \ee
So by (\ref{y54}), $\rho$ is a (closed and ) conformal 1-form with
respect to $h$. Therefore, by Theorem \ref{th1}, we get the proof
of Theorem \ref{th61} in this case.

\

 \noindent {\bf Case II:}
 Suppose that (\ref{y53}) holds. In this case, let $u,v,w$  in (\ref{y1})
be defined by
 $$
 u(b^2)=\frac{(1\mp b^2)^3}{(1\pm 2b^2)^{3/2}}, \ \
 v(b^2)=\frac{9}{8b^2}\Big\{(1\pm 2b^2)^{3/2}-\frac{1\mp
 2b^2+4b^4}{(1\pm 2b^2)^{3/2}}\Big\},\ \ w(b^2)=1.
 $$
Then it is shown in \cite{Y3}  that (\ref{y53}) is equivalent to
 \be\label{y55}
 \widetilde{r}_{ij}=\frac{\theta (1\mp b^2)^2}{(1\pm
 2b^2)^{5/2}}\ h_{ij},
 \ee
 where we define $\widetilde{r}_{ij}:=(p_{i|j}+p_{j|i})/2$ for the
covariant derivative of $\rho$ with respect to $h$. So by
(\ref{y55}), $\rho$ is a  conformal 1-form with respect to $h$.
Therefore, by Theorem \ref{th1}, we get the proof of Theorem
\ref{th61} in this case.  \qed

\section{Projectively flat metrics}\label{sec7}

For  a locally projectively flat $(\alpha,\beta)$-metric
$F=\alpha\phi(\beta/\alpha)$ of non-Randers type with the
dimension $n\ge 3$, the local solution of a conformal vector field
has been obtained in \cite{Y8} (also see Corollary \ref{cor41}).

In the following, we consider the conformal vector fields on
locally projectively flat Randers spaces. In this case, we cannot
determine the local structure of those conformal vector fields,
and the conformal vector fields on such spaces are not necessarily
homothetic.

A Randers metric $F=\alpha+\beta$ is locally projectively flat if
and only if $\alpha$ is of constant sectional curvature and
$\beta$ is closed (\cite{BaMa}). In this case, we don't need to
deform $\alpha$ and $\beta$ by (\ref{y1}), or namely, we just
choose $u=w=1,v=0$ in (\ref{y1}). To study conformal vector fields
on locally projectively flat Randers spaces, we first show a
lemma.

\begin{lem}\label{lem71}
 Let $F=\alpha \phi(s)$, $s=\beta/\alpha$, be an $n(\ge
2)$-dimensional non-Riemannian
  $(\alpha,\beta)$-metric on the manifold $M$, where $\phi(0)=1$.
  Suppose $\beta$ is closed and $V$ is a conformal vector field
  with the conformal factor $c$. Then we have $c_i=\theta b_i$,
  where $c_i:=c_{x^i}$ and $\theta=\theta(x)$ is a scalar
  function.
\end{lem}

{\it Proof :} By Lemma \ref{lem22}, we have (\ref{y13}). Then
since $\beta$ is closed, the second equation in (\ref{y13}) is
written as
 $$
 H_i=-2cb_i,\ \ \ (H:=b^iV_i,\ \ H_i:=H_{x^i}).
 $$
Therefore, we have
 $$0=H_{i;j}-H_{j;i}=(-2c_jb_i-2cb_{i;j})-(-2c_ib_j-2cb_{j;i})=-2(b_ic_j-b_jc_i),$$
 which imply $c_ib_j=c_jb_i$. Thus we have $c_i=\theta b_i$ for
 some scalar function $\theta=\theta(x)$ since $\beta\ne0$.  \qed

\

Now let $F=\alpha+\beta$ be  locally projectively flat on a
manifold $M$ with the dimension $n\ge 2$. Locally express
 $$
 \alpha=\frac{2}{1+\mu|x|^2}|y|,
 $$
  and put $V=V^i\pa/\pa x^i$ by (\ref{y16}), namely,
 \be\label{y56}
  V^i=-2\big(\tau+\langle
 \eta,x\rangle\big)x^i+|x|^2\eta^i+q_r^ix^r+\gamma^i.
 \ee
By Lemma \ref{lem21}, $V$ is a conformal vector field of
$(M,\alpha)$ with the conformal factor $c$ given by (\ref{y17}),
namely,
 \be\label{y57}
 c=\frac{\tau(1-\mu |x|^2)+ \langle \mu
\gamma+\eta,x\rangle}{1+\mu |x|^2}
 \ee
By Lemma \ref{lem22}, if the above $V$ is also a conformal vector
field of $(M,F)$, then we must have the second equation in
(\ref{y13}), namely,
 \be\label{y58}
 V^j\frac{\pa b_i}{\pa
 x^j}+b_j\frac{\pa V^j}{\pa x^i}=-  2cb_i.
 \ee
By Lemma \ref{lem71}, we may put
 \be\label{y59}
 b_i=f(c)c_i,
 \ee
where $f$ is some function. Now plug (\ref{y56}), (\ref{y57}) and
(\ref{y59}) into (\ref{y58}), and then we obtain
 \be\label{y60}
 A_1|x|^2+A_2=0,
 \ee
where
 \beqn
 A_1\hspace{-0.5cm}&&=\big(2\mu(c+\tau)x^i-\eta^i-\mu\gamma^i\big)(2\mu
 c^2-2\mu\tau
 c-4\mu\tau^2-|\eta|^2-\mu\langle\eta,\gamma\rangle)f'(c)\\
 &&+\mu\big[2(\tau-c)\eta^i-2\mu(\tau+c)\gamma^i+4\mu(c^2-\tau^2)x^i-q^i_k\eta^k-\mu
 q^i_k\gamma^k\big]f(c),\\
 A_2\hspace{-0.5cm}&&=\Big\{2\mu\big[\tau
 f'(c)+cf'(c)+f(c)\big]\langle 4\mu\tau\gamma+Q\eta+\mu
 Q\gamma,x\rangle-2\mu(\tau+c)(2\tau c-2c^2\\
 &&+\mu|\gamma|^2+\langle
 \eta,\gamma\rangle)f'(c)-2(\mu^2|\gamma|^2-|\eta|^2-2\mu\tau^2-2\mu
 c^2)f(c)\Big\}x^i-\\
 &&(\eta^i+\mu\gamma^i)f'(c)\langle 4\mu\tau\gamma+Q\eta+\mu
 Q\gamma,x\rangle+(\eta^i+\mu\gamma^i)(2\tau c-2c^2
 +\mu|\gamma|^2+\langle
 \eta,\gamma\rangle)f'(c)\\
 &&+\big[2(\tau-c)\eta^i-2\mu(\tau+c)\gamma^i-q^i_k\eta^k-\mu
 q^i_k\gamma^k\big]f(c),
 \eeqn
To construct non-homothetic vector field, we assume $A_1=A_2=0$.
Then by the expression of $A_2$, we first assume
 \be\label{y61}
 Q\eta=-\mu(4\tau\gamma+Q\gamma).
 \ee
Plug (\ref{y61}) into $A_2$, and we obtain
 \beq
 0\hspace{-0.5cm}&&=-2\big[\mu(\tau+c)(2\tau c-2c^2+\mu|\gamma|^2+\langle
 \eta,\gamma\rangle)f'(c)+(\mu^2|\gamma|^2-|\eta|^2-2\mu\tau^2-2\mu
 c^2)f(c)\big]x^i\nonumber\\
 &&(\eta^i+\mu\gamma^i)\big[(2\tau c-2c^2+\mu|\gamma|^2+\langle
 \eta,\gamma\rangle)f'(c)-2(c-\tau)f(c)\big].\label{y62}
 \eeq
By the second term of the right hand side of (\ref{y62}), we may
assume
 \be\label{y63}
 f'(c)=\frac{2(c-\tau)}{2\tau c-2c^2+\mu|\gamma|^2+\langle
 \eta,\gamma\rangle}f(c).
 \ee
Then we can obtain the function $f(c)$ by solving the ODE
(\ref{y63}). Further, plugging (\ref{y63}) into (\ref{y62}) gives
  \be\label{y64}
 |\eta|^2=\mu(\mu|\gamma|^2-4\tau^2).
  \ee
Conversely, if we assume (\ref{y61}), (\ref{y63}) and (\ref{y64})
hold, we can directly verify that $A_1=A_2=0$.

Now by the above construction, we obtain the following example.

\begin{ex}\label{ex72}
 Let $F=\alpha+\beta$ be a Randers metric on a manifold $M$. Locally, define
  \beqn
 \alpha:&&\hspace{-0.6cm}=\frac{2}{1+\mu|x|^2}|y|, \ \ \ \ \ \
 \beta:=f(c)c_{x^i},\\
  c:&&\hspace{-0.6cm}=\frac{\tau(1-\mu |x|^2)+ \langle \mu
\gamma+\eta,x\rangle}{1+\mu |x|^2},\\
 V^i:&&\hspace{-0.6cm}=-2\big(\tau+\langle
 \eta,x\rangle\big)x^i+|x|^2\eta^i+q_r^ix^r+\gamma^i,
  \eeqn
  where $f$ is a function,  $\mu\ (\ne0)$ and $\tau$ and  $\eta=(\eta^i)$ and
  $\gamma=(\gamma^i)$ are of constant values, and $Q=(q^i_j)$ is a
  constant skew-symmetric matrix, and those  parameters
  satisfy (\ref{y61}), (\ref{y63}) and (\ref{y64}). It is easy to
  see that $F$ is locally projectively flat.

 The above construction has shown that the vector field $V=(V^i)$ is conformal on $(M,F)$
   with the conformal factor $c$, and $V$ is not homothetc.
\end{ex}

In Example \ref{ex72}, if we take $Q=0$ and $\gamma=0$, then we
have
 \beqn
   |\eta|^2\hspace{-0.5cm}&&=-4\mu\tau^2, \ \ \ \ \  c=\frac{\tau(1-\mu |x|^2)+ \langle \eta,x\rangle}{1+\mu |x|^2},\\
 \beta\hspace{-0.5cm}&&=\frac{1}{\tau(1-\mu|x|^2)+\langle
 \eta,x\rangle}\Big\{\langle \eta,y\rangle-\frac{2\mu(2\tau+\langle \eta,x\rangle)\langle x,
 y\rangle}{1+\mu|x|^2}\Big\},\\
 V^i\hspace{-0.6cm}&&=-2\big(\tau+\langle
 \eta,x\rangle\big)x^i+|x|^2\eta^i,
 \eeqn
where we have put $f(c)=1/c$ (by (\ref{y63})). Thus we obtain a
simple family of non-homothetic conformal vector fields on a
corresponding family of locally projectively flat Randers spaces,
which have been shown in \cite{Y8}.

In a more general case, Example \ref{ex72} also shows the
following remark.

\begin{rem}\label{rem73}
 In Theorem \ref{th1}, if we assume $h$ is of constant non-zero
 sectional curvature, and $\rho$ is just closed (instead of being
 conformal), then the conformal vector field of $(M,F)$ is not
 necessary homothetic.
\end{rem}

\vspace{0.6cm}

\noindent Guojun Yang \\
Department of Mathematics \\
Sichuan University \\
Chengdu 610064, P. R. China \\
yangguojun@scu.edu.cn

\end{document}